\documentclass{amsart}
\usepackage{amsmath}
\usepackage{amscd}
\usepackage{amssymb}

\renewcommand{\div}{\operatorname{div}}

\newtheorem{theorem}{Theorem}[section]
\newtheorem{theorem/definition}{Theorem/Definition}[section]
\newtheorem{proposition}{Proposition}[section]
\newtheorem{lemma}{Lemma}[section]
\newtheorem{corollary}{Corollary}[section]
\theoremstyle{remark}
\newtheorem{remark}{Remark}[section]
\theoremstyle{definition}

\makeatletter
\@namedef{subjclassname@2020}{%
  \textup{2020} Mathematics Subject Classification}
\makeatother

\begin{document}
\title
{On Complete Gradient Steady Ricci Solitons with vanishing $D$-tensor}
\author{Huai-Dong Cao$^{\dag}$ and Jiangtao YU}
\address{Department of Mathematics\\ Lehigh University\\
Bethlehem, PA 18015} \email{huc2@lehigh.edu}

\address{Department of Mathematics\\ Lehigh University\\
Bethlehem, PA 18015} \email{jiy314@lehigh.edu}

\thanks{$^{\dag }$Research partially supported  by Simons Foundation Collaboration Grant \#586694 HC}

\begin{abstract}
In this paper, we extend the work in \cite{CC13} to classify $n$-dimensional ($n\ge 5$) complete
$D$-flat gradient steady Ricci solitons. More precisely, we
prove that any $n$-dimensional complete noncompact gradient steady Ricci
soliton with vanishing D-tensor is either Ricci-flat, or isometric to the Bryant soliton.  
Furthermore, the proof extends to  the shrinking case and the expanding case as well.

\end{abstract}

\maketitle
\date{}

\section{Introduction}

A complete Riemannian manifold $(M^n,  g_{ij})$ is called a {\it gradient steady Ricci soliton} if there exists
a smooth function $F$ on $M^n$ such that the Ricci tensor $R_{ij}$
of the metric $g_{ij}$ is given by the Hessian of $F$:
$$R_{ij}=\nabla_i\nabla_jF. \eqno(1.1)$$
Such a function $F$ is called a {\it potential function} of the gradient steady soliton.
Clearly, when $F$ is a constant the gradient steady Ricci soliton is simply a Ricci flat manifold. 
Thus Ricci solitons are natural extensions of Einstein metrics. 
Gradient steady solitons play an important role in
Hamilton's Ricci flow as they correspond to translating
solutions and often arise as Type II singularity models. Therefore, one is interested
in classifying them and understanding their geometry.

It turns out that compact steady Ricci solitons must be Ricci flat. In dimension $n=2$, Hamilton \cite{Ha88} discovered the 
first example of a complete noncompact gradient steady
soliton on $\mathbb R^2$, called the {\it cigar soliton}, where the metric is given explicitly by
$$ ds^2=\frac{dx^2 +dy^2}{1+x^2+y^2},$$
and the potential function is $F=\log (1+x^2+y^2)$.  
The cigar soliton has positive curvature and is asymptotic to a cylinder of finite
circumference at infinity.  Furthermore, Hamilton \cite{Ha88} showed that the cigar soliton is the only
complete steady soliton on a two-dimensional manifold with
bounded (scalar) curvature $R$ which assumes its maximum
at the origin. For $n\geq 3$,  Bryant proved in late 1980s that
there exists, up to scalings, a unique complete rotationally symmetric gradient steady Ricci
soliton on $\Bbb R^n$; see \cite{Bryant} and Chow et al \cite{Chow et al 1}
for a detailed description. The Bryant soliton has positive sectional curvature, linear curvature decay
and volume growth on the order of $r^{(n+1)/2}$. Here, $r$ denotes the geodesic distance from the origin. 
In the K\"ahler case,
the first author \cite{Cao94} constructed a complete $U(m)$-invariant gradient steady K\"ahler-Ricci soliton on $\mathbb{C}^m$, for $m\geq 2$, with positive sectional curvature, linear curvature decay, and volume growth on the order of $r^m$.

In dimension $n=3$, Brendle \cite {Brendle13} confirmed a conjecture of Perelman that the Bryant soliton is the
only complete noncompact $\kappa$-noncollapsed gradient steady
soliton with positive sectional curvature. Subsequently,  he 
extended his result to dimension $n\ge 4$ under the extra assumption that the steady soliton is asymptotically cylindrical \cite {Brendle14}; see also related recent works by Deng and Zhu \cite{DZ19, DZ20}. Naturally,  for $n=3$, the remaining open problem is to understand complete noncompact collapsed gradient steady Ricci solitons with positive curvature. We refer the reader to a recent paper \cite{CHe18} for more information on this 
subject\footnote{Note added in proof: Very recently Y. Lai \cite{Lai} has proved the existence of the flying wing $3$-dimensional gradient steady Ricci soliton.}.

On the other hand, for $n\ge 3$, Cao and Chen \cite{CC09} proved that the Bryant
soliton is the only complete noncompact positively curved
{\it locally conformally flat} gradient steady soliton; see also the work of Catino and Mantegazza \cite{CM11} for an independent proof when $n\ge 4$. Furthermore, in dimension $n=4$,  Chen and Wang \cite{CW} showed that half-conformal flatness implies locally conformal flatness, thus improving the result of \cite{CC09, CM11} for $n=4$. 

In \cite{CC09, CC13}, a key covariant 3-tensor $D$ was introduced to study locally conformally flat and Bach-flat gradient Ricci solitons. 
Specifically, the $D$-tensor is defined by
$$ D_{ijk} = \frac{1} {n-2} (A_{ik} \nabla_jF-A_{jk} \nabla_i F) + \frac {1} {(n-1)(n-2)} (g_{ik}E_{jl} -g_{jk}E_{il})\nabla_lF, \eqno(1.2)$$
where $A_{ij}$ is the Schouten tensor and $E_{ij}$ is the Einstein tensor.  Equivalently, as shown in \cite{CC13}, 
$$D_{ijk}=C_{ijk}-W_{ijkl}\nabla_l F, \eqno(1.3)$$
where $C_{ijk}$ denotes the Cotton tensor and $W_{ijkl}$ is the Weyl tensor. Clearly, locally conformal flatness ($W_{ijkl}=0$) implies $D_{ijk}=0$. In \cite{CC13}, it was shown that $D_{ijk}=0$ implies $W_{ijkl}=0$ in dimension $n=4$ for gradient Ricci solitons. Consequently, any $4$-dimensional complete noncompact gradient steady Ricci soliton with vanishing D-tensor is either flat or isometric to the Bryant soliton. In addition, they provided several equivalent conditions characterizing $D_{ijk}=0$ in dimension $n\ge 5$; see Proposition 2.3.  However, it remained unclear what is the structure of complete noncompact gradient steady Ricci solitons with vanishing $D$-tensor for $n\ge 5$. 

In this paper, we extend the work of Cao-Chen \cite{CC13} to classify $n$-dimensional complete noncompact $D$-flat gradient steady Ricci solitons for $n\ge 5$.

\begin{theorem} Let $(M^n, g_{ij}, F)$, $n\ge 5$, be a complete noncompact gradient steady Ricci soliton with vanishing $D$-tensor. Then $(M^n, g_{ij}, F)$ is either Ricci-flat with a constant potential function, or a quotient of the product steady soliton $N^{n-1} \times {\mathbb R}$, where $N^{n-1}$ is Ricci-flat, or isometric to the Bryant soliton (up to scalings).
\end{theorem}

\medskip
\noindent 
{\bf Corollary 1.2.} {\em Let $(M^n, g_{ij}, F)$, $n\ge 5$, be a complete noncompact 
gradient steady Ricci soliton with positive scalar curvature and vanishing $D$-tensor. Then $(M^n, g_{ij}, F)$ is  isometric to the Bryant soliton.}  

\begin{remark} Note that Theorem 1.1 does not require any a priori curvature sign or bound conditions. Also, the proof can be adapted to treat gradient shrinking or expanding Ricci solitons with vanishing $D$-tensor; see Propositions  3.2-3.3. 
\end{remark}

\begin{remark}
Recently, Kim \cite{Kim} showed  that any $4$-dimensional complete noncompact gradient steady soliton with harmonic Weyl curvature (i.e., $C_{ijk}=0$) is either Ricci-flat or isometric to the Bryant soliton. In a very recent preprint, Fengjiang Li \cite{Li} has extended Kim's result to all dimensions $n\ge 5$ where our Theorem 1.1 is used to treat one of the cases in the classification. 
\end{remark}

\begin{remark} Corollary 1.2 was first proved in \cite{Cao et al} under the stronger assumption that the Ricci curvature is positive and the scalar curvature attains its maximum at some interior point. 

\end{remark}

\medskip \noindent {\bf Acknowledgements.} We would like to thank Chenxu He and Junming Xie for helpful discussions.

\section{Preliminaries}

In this section, we fix our notations and recall some basic facts
and known results about gradient Ricci solitons that we shall need later.

\subsection{Classical conformal tensors} First of all, we recall that on any $n$-dimensional Riemannian
manifold $(M^n, g_{ij} )$ ($n\ge 3$),  the Weyl curvature tensor
is given by

\begin{align*}
W_{ijkl}  = & R_{ijkl} - \frac{1}{n-2}(g_{ik}R_{jl}-g_{il}R_{jk}-g_{jk}R_{il}+g_{jl}R_{ik})\\
& + \frac{R}{(n-1)(n-2)} (g_{ik}g_{jl}-g_{il}g_{jk}), \\
\end{align*}
and  the Cotton tensor by
 $$C_{ijk}=\nabla_i R_{jk}-\nabla_j R_{ik}-\frac {1}{2(n-1)} ( g_{jk} \nabla_i R - g_{ik} \nabla_j R).$$
 
\begin{remark} In terms of the Schouten tensor 
$$ A_{ij}= R_{ij}- \frac{R}{2(n-1)}g_{ij}, \eqno(2.1)$$  we have  
$$ W_{ijkl}=R_{ijkl} - \frac{1}{n-2} (g_{ik}A_{jl}-g_{il}A_{jk}-g_{jk}A_{il}+g_{jl}A_{ik}), $$ and
$$C_{ijk}=\nabla_i A_{jk}-\nabla_j A_{ik}.$$
\end{remark}
 
 It is well known that, for $n=3$, $W_{ijkl}$ vanishes identically, while $C_{ijk}=0$ if and only if
$(M^3, g_{ij})$ is locally conformally flat; for $n\ge 4$,
$W_{ijkl}=0$ if and only if  $(M^n, g_{ij})$ is locally
conformally flat.  Moreover, for $n\ge 4$, the Cotton tensor
$C_{ijk}$ is, up to a constant factor, the divergence of the Weyl
tensor:
$$C_{ijk}=-\frac{n-2}{n-3} \nabla_l W_{ijkl}, \eqno(2.2)$$
hence the vanishing of the Cotton tensor $C_{ij k} = 0$ (in
dimension $n\ge 4$) is also referred to as being harmonic Weyl.

Note that  $C_{ijk}$ is skew-symmetric in the
first two indices and trace-free in any two indices:
$$ C_{ijk}=-C_{jik} \quad \mbox{and} \quad g^{ij}C_{ijk}=g^{ik}C_{ijk}=0.$$

Moreover, for $n\ge 4$, the Bach tensor is defined by
$$B_{ij}=\frac 1 {n-3}\nabla^k\nabla^l W_{ikjl}+\frac 1 {n-2}
R_{kl}W{_i}{^k}{_j}^l. \eqno(2.3)$$ By (2.2), we have
$$B_{ij}=\frac 1 {n-2} (\nabla_k C_{kij} + R_{kl}W{_i}{^k}{_j}^l). \eqno(2.4)$$

\subsection{Two basic facts for gradient steady Ricci solitons}  Next, we recall two  basic facts about  gradient steady Ricci solitons.

\begin{lemma} {\bf (Hamilton \cite{Ha95F})} Let $(M^n, g_{ij}, F)$
be a  gradient steady Ricci soliton satisfying Eq. (1.1). 
Then, we have
$$\nabla_iR=-2R_{ij}\nabla_jF, \eqno(2.5)$$ and
$$R+|\nabla F|^2=C_0 \eqno(2.6)$$ for some constant $C_0$. Here $R$
denotes the scalar curvature.
\end{lemma}

\begin{lemma} {\bf (B.-L. Chen \cite{BChen})} Let $(M^n, g_{ij}, F)$ be a complete gradient steady soliton. Then it has nonnegative
scalar curvature $R\ge 0$.
\end{lemma}

Lemma 2.2 is a special case of a more general result of B.-L. Chen
\cite{BChen} which states that $R\ge 0$ for any complete ancient solution
to the Ricci flow.

\subsection {The covariant 3-tensor $D_{ijk}$}

For any gradient Ricci soliton  satisfying the defining equation 
$$R_{ij}+\nabla_i\nabla_jf=\rho g_{ij}, \eqno (2.7)$$ with any $\rho \in {\mathbb R}$, Cao-Chen \cite{CC09, CC13} introduced a covariant 3-tensor $D_{ijk}$ defined by
\begin{align*}
 D_{ijk} = & \frac{1}{n-2}(R_{jk} \nabla_i f- R_{ik} \nabla_j f) +\frac{1}{2(n-1)(n-2)} (g_{jk}\nabla_i R-g_{ik} \nabla_j R)\\
              & - \frac{R}{(n-1)(n-2)} (g_{jk}\nabla_i f  - g_{ik} \nabla_j f ).
 \end{align*}
Note that, by using (2.5) which is valid for (2.7) with $f=-F$, $D_{ijk}$ can also be expressed as 
$$ D_{ijk} = \frac{1} {n-2} (A_{jk} \nabla_i f- A_{ik} \nabla_jf) + \frac {1} {(n-1)(n-2)} (g_{jk}E_{il} - g_{ik}E_{jl})\nabla_lf, \eqno(2.8)$$
where $A_{ij}$ is the Schouten tensor in (2.1) and $E_{ij}=R_{ij}-\frac {R} {2} g_{ij}$ is the Einstein tensor. 

This 3-tensor $D_{ijk}$ is closely tied to the Cotton tensor, as well as the Bach tensor, and 
played a significant role on classifying
locally conformally flat gradient steady solitons in  \cite{CC09} and Bach flat shrinking Ricci solitons in \cite{CC13}.

\begin{lemma} {\bf (Cao-Chen \cite{CC13})} 
Let $(M^n, g_{ij}, f)$ ($n\ge 3$) be a complete gradient soliton
satisfying (2.7). Then $D_{ijk}$ is related to the Cotton tensor
$C_{ijk}$ and the Weyl tensor $W_{ijkl}$ by
$$ D_{ijk}=C_{ijk}+W_{ijkl}\nabla_l f.$$
\end{lemma}

\begin{remark}
By Lemma 2.3, it follows that $D_{ijk}$ is equal to the Cotton tensor $C_{ijk}$ in dimension $n=3$.  In addition, for $n\ge 3$, it is easy to see that 
$$D_{ijk}\nabla_k f=C_{ijk} \nabla_k f.$$
Also, $D_{ijk}$ vanishes if $(M^n, g_{ij}, f)$ ($n\ge 3$) is either Einstein with a constant potential function $f$, or locally conformally flat. 
Moreover, like the Cotton tensor $C_{ijk}$, $D_{ijk}$ is skew-symmetric in the first two indices and trace-free in any 
two indices:   
$$ D_{ijk}=-D_{jik} \quad \mbox{and} \quad g^{ij}D_{ijk}=g^{ik}D_{ijk}=0. $$ 
\end{remark}

What is so special about $D_{ijk}$ is the following key identity,
which links  the norm of $D_{ijk}$ to the geometry of the level
surfaces of the potential function $f$. 

\begin{proposition} {\bf (Cao-Chen \cite{CC09, CC13})}  Let $(M^n, g_{ij}, f)$ ($n\ge 3$) be an $n$-dimensional
gradient Ricci soliton satisfying (2.7). Then, at any point $p\in
M^n$ where $\nabla f(p)\ne 0$, we have
$$|D_{ijk}|^2=\frac{2|\nabla f|^4}{(n-2)^2} |h_{ab}-\frac{H}{n-1}g_{ab}|^2 +\frac{1}{2(n-1)(n-2)}
|\nabla_a R|^2,$$ where $h_{ab}$ and $H$ are the second
fundamental form and the mean curvature of the level surface
$\Sigma=\{f=f(p)\}$ respectively, and $g_{ab}$ is the induced metric on $\Sigma$.
\end{proposition}

\begin{remark} For any gradient Ricci solitons satisfying (2.7), the Bach tensor is related to the D-tensor and the Cotton tensor as follows. 
$$B_{ij} =-\frac 1{n-2}(\nabla_k D_{ikj}+\frac{n-3}{n-2}C_{jli}\nabla_l f).$$
In particular, by  \cite{CC13}, the vanishing of $D$-tensor implies the vanishing of the Bach tensor $B$. Conversely, it was shown in \cite{CC13} (or \cite{Cao et al}) that Bach flat gradient shrinking Ricci solitons (or Bach-flat gradient steady solitons with positive Ricci curvature) must be $D$-flat. More generally, we have the following characterizations of $D=0$. 
\end{remark}

\begin{proposition} {\bf (Cao-Chen \cite{CC13})} Let $(M^4, g_{ij}, f)$ be a complete non-trivial gradient Ricci soliton satisfying (2.7) and with $D_{ijk}=0$.  Then  $(M^4, g_{ij}, f)$ is locally conformally flat, i.e., $W_{ijkl}=0$. 
\end{proposition}

\begin{proposition} {\bf (Cao-Chen \cite{CC13})} Let $(M^n, g_{ij}, f)$ ($n\ge 5$) be a nontrivial gradient Ricci soliton satisfying (2.7). Then the following statements are equivalent: 

\smallskip
(a) $D_{ijk}=0$;

\smallskip
(b) $C_{ijk}=0$ and $W_{1ijk}=0$ for $1\le i,j,k\le n$; 

\smallskip
(c) $C_{ijk}\nabla_i f=0$ and $W_{1a1b}=0$  for $1\le i,j,k\le n$ and $2\le a, b\le n$.

\smallskip
(d) $\div B\cdot \nabla f=0$ and $W_{1a1b}=0$  for $2\le a, b\le n$.
\end{proposition}  

\begin{remark} Note that in general $(b)\Longrightarrow (c) \Longrightarrow (d)$ for any Riemannian manifold.   
\end{remark}

Combining Proposition 2.2 and the classification for locally conformally flat gradient steady Ricci soliton \cite{CC09, CM11}, we have 

\begin{corollary} {\bf (Cao-Chen \cite{CC13})} Let $(M^4, g_{ij}, f)$  be a complete gradient {\sl steady} Ricci soliton with $D_{ijk}=0$,  then $(M^4, g_{ij}, f)$ is either flat or isometric to the Bryant soliton.
\end{corollary}

\section{The proof of Theorem 1.1}

Throughout this section,  we assume that $(M^n, g_{ij}, F)$ ($n\ge
4$) is a complete gradient steady soliton satisfying (1.1). 

First of all, we shall need the following result from \cite{CC13}  (and \cite{Cao et al}) which shows that the vanishing of $D_{ijk}$ implies
many nice properties about the geometry of a Ricci soliton and the regular level surfaces of its potential
function.

\begin{proposition}  {\bf (Cao-Chen \cite{CC13})} 
Let $(M^n, g_{ij}, f)$ ($n\ge 4 $) be any complete gradient Ricci
soliton satisfying (2.7) with  $D_{ijk}=0$.  Let $c$ be a regular value of the potential function $f$
and $\Sigma_c=\{f=c\}$ be the level surface of $f$. Set
$e_1=\nabla f /|\nabla f |$ and pick any orthonormal frame $e_2,
\cdots, e_n$ tangent to the level surface $\Sigma_c$. Then

 \medskip

(a) $|\nabla f|^2$ and the scalar curvature $R$ of $(M^n,
g_{ij}, f)$ are constant on each connected component of $\Sigma_c$;

\smallskip
(b) $R_{1a}=0$ for $2\le a\le n$, and $e_1=\nabla f /|\nabla f |$
is an eigenvector of the Ricci tensor $Rc$;

\smallskip

(c) the second fundamental form $h_{ab}$ of $\Sigma_c$ is of the
form $h_{ab}=\frac{H}{n-1} g_{ab}$; 

\smallskip

(d) the mean curvature $H$ is constant on each connected component of $\Sigma_c$;

\smallskip

(e) on $\Sigma_c$, the Ricci tensor of $(M^n, g_{ij}, f)$ either
has a unique eigenvalue $\lambda$, or has two distinct eigenvalues
$\lambda$ and $\mu$ of multiplicity $1$ and $n-1$ respectively. In
either case, $e_1=\nabla f /|\nabla f |$ is an eigenvector of
$\lambda$.

\smallskip
(f) locally, on some neighborhood $U$ of the level surface $\Sigma_c$ where $\nabla f\neq 0$, $g_{ij}$ is a warped product metric of the form $$ds^2=dr^2+\varphi^2(r)\bar g_{\Sigma_c},$$ where $\bar g_{\Sigma_c}$ denotes the induced metric on $\Sigma_c$. Furthermore, $(\Sigma_c, \bar g_{\Sigma_c})$ is Einstein. 
\end{proposition}

\bigskip

\noindent {\bf Proof of Theorem 1.1}:  Let $(M^n,
g_{ij}, F)$, $n\ge 4$, be a complete noncompact gradient steady Ricci soliton.  Without loss of generality, we  assume $M^n$ is simply connected, otherwise we can consider its universal cover $\tilde M^n$ and lift $g_{ij}$ and $F$ to $\tilde M^n$. We divide the proof into two cases: 

\medskip
{\bf Case 1:}  $|\nabla F|^2=0$ on some nonempty open set of $M^n$. 

\smallskip
In this case, since any gradient Ricci soliton is analytic
in harmonic coordinates (see, e.g., \cite{Ivey}), it follows that $|\nabla F|^2=0$ on $M$ everywhere,
i.e., $F$ is a constant function and $(M^n, g_{ij})$ is Ricci-flat.

\medskip
{\bf Case 2:}  The open set $\Omega =\{p\in M | \nabla F(p)\ne 0\}$ is dense. 

\smallskip

Pick a regular value $c_0$ of the potential function $F$ and consider the level surface $\Sigma_{c_0}=F^{-1}(c_0)$. Suppose $I\subset {\mathbb R}$ is an open interval containing $c_0$ such that
$F$ has no critical points in the open neighborhood $U_{I}=F^{-1} (I)$ of $\Sigma_{c_{0}}$.
In this case,  by applying Proposition 3.1 with $F=-f$ and $\rho=0$, one obtains a local warped product structure on $U_I$ as mentioned in \cite {CC13} (Remark 3.2) and detailed in \cite{Cao et al}:  
$$ ds^2=dr^2+\varphi^2(r)\bar g_{\Sigma_{c_0}}, \eqno(3.1)$$ where  
$$ r(x)=\int_{F=c_0}^F \frac{dF} {|\nabla F|}, \eqno (3.2)
$$ 
and $ {\bar g}:= {\bar g}_{\Sigma_{c_{0}}}=g_{ab}(r, \theta) d\theta^a d\theta^{b} $
is the induced metric on $\Sigma_{c_{0}}$, for any local coordinates system $\theta=(\theta^2, \cdots, \theta^n)$ of $\Sigma_{c_{0}}$. Furthermore, $(\Sigma_{c_0},  \bar g)$ is Einstein.  

Let $\nabla r=\frac{\partial} {\partial r}$, then we have $ |\nabla r|=1$ and $\nabla F=F'(r) \frac{\partial} {\partial r}$ on $U_{I}$. Note that  $F'(r)$ does not change sign on $U_{I}$ because  $F$ has no critical points there. Thus we may assume $I=\{\alpha<r<\beta\}=(\alpha, \beta)$
with  $F'(r)>0$ on $I$.  It is also easy to check that
$$ \nabla_{\frac {\partial} {\partial r}} \frac{\partial} {\partial r}=0, \eqno (3.3)
$$
so integral curves to $\nabla r$ are normal geodesics. Since we assume $M^n$ is simply connected, $\Sigma_{c_0}=r^{-1}(r_0)$ is a two-sided embedded hypersurface and $r$ is locally a (signed) distance function from $\Sigma_{c_0}$. 

\medskip
{\bf Claim 1.} $(M^n, g_{ij})$ is isometric to one of the following warped products: 

\medskip
(i) \ $({\mathbb R}, \,dr^2)\times\!{_\varphi}(\Sigma_{c_0}, \, \bar g) $; or

\smallskip
(ii) $([0, \infty), \,dr^2)\times\!{_\varphi}(\Sigma_{c_0}, \, \bar g). $
\medskip

Set $I_{\max}\supseteq I$ to be the maximal interval for which the warped product structure (3.1) holds. Now, suppose $|\nabla F|(\bar r)=0$ for some $\bar r$, with either $(\bar r-\epsilon, \bar r) \subset I_{\max}$ or $(\bar r, \bar r+\epsilon) \subset I_{\max}$. As pointed out in \cite{FG2}, by continuity and smoothness, as long as the warping function $\varphi (\bar r)\ne 0$, $\Sigma_{\bar c}=F^{-1}(\bar c)=r^{-1}(\bar r)$ is a totally umbilical submanifold of $(M^n, g_{ij})$ and the warped structure (3.1) can be extended to $(\bar r-\epsilon, \bar r+\epsilon) \subset I_{\max}$ for some $\epsilon >0$ sufficient small.  Here, we have used the assumption that the set of critical points of $F$ is nowhere dense in Case 2. 
Now it follows that either $\varphi$ never vanishes so that $I_{\max}=(-\infty, \infty)$, or it vanishes at exactly one value of $r$ so that, e.g, $I_{\max}=[a, \infty)$,  or at exactly two values of $r$ so that $I_{\max}=[a, b]$. However, the last case does not happen since $M^n$ is noncompact by assumption. Obviously, by a suitable shifting in $r$, we can assume $[a, \infty)=[0, \infty)$. This proves Claim 1.

\medskip
{\bf Claim 2.} If $(M^n, g_{ij})=([0, \infty), \,dr^2)\times\!{_\varphi}(\Sigma_{c_0}, \, \bar g), $ then $(\Sigma_{c_0}, \, \bar g)$ is isometric to a round sphere $({\mathbb S}^{n-1}, \, \bar g_0)$ and $(M^n, g_{ij}, F)$ is isometric to the Bryant soliton. 

\medskip
This follows from the fact that $r(x)$ is simply the distance function
$d (x_0,x)$ from $x_0=\varphi^{-1}(0)$. So level surfaces of $F$ are geodesic spheres centered at $x_0$ which are diffeomorphic to the $(n -1)$-sphere ${\mathbb S}^{n-1}$. In addition, either by an argument in \cite{Cao et al} or by the smoothness of the metric g at $x_0$ (see, e.g., Lemma 9.114 in \cite{Besse}), the induced Einstein metric $\bar g$ on $\Sigma_{c_0}$ is necessary round.  As a consequence, $(M^n, g_{ij})$ is rotationally symmetric, hence isometric to the Bryant soliton by \cite{CC09, CM11}.

\medskip
{\bf Claim 3.} If $(M^n, g_{ij})=({\mathbb R}, \,dr^2)\times\!{_\varphi}(\Sigma_{c_0}, \, \bar g) $, then 
$\varphi $ is a constant function and  
$$(M^n, g)=({\mathbb R}, \,dr^2)\times\,(\Sigma_{c_0}, \, \bar g), $$ where $(\Sigma_{c_0}, \, \bar g)$ is  necessarily Ricci-flat. 

\medskip
First of all, from either direct computations or \cite{Besse, Oneill}, the Ricci tensor of $(M^n, g)=({\mathbb R}, \,dr^2)\times\!{_\varphi}(\Sigma_{c_0}, \, \bar g) $ is given by
$$
R_{11}=-(n-1)\frac {\varphi''} {\varphi}, \quad R_{1a}=0 \quad (2\le a\le n), \eqno(3.4)
$$ 
$$
R_{ab}=\bar R_{ab}-\big [\varphi\varphi'' +(n-2) (\varphi')^2\big] \bar g_{ab}  \quad (2\le a,b\le n), \eqno(3.5)
$$
where $\bar R_{ab}$ denotes the Ricci tensor of $ \bar g_{ab}$. Moreover, the Hessian of $F$ can be computed by 
$$ \nabla_1\nabla_1 F=F'', \quad \nabla_a\nabla_1 F=0, \quad \nabla_a\nabla_b F=\varphi\varphi'F' \bar g_{ab}.\eqno (3.6)$$
It follows that the scalar curvature $R$ of $(M^n, g)$ is given by
$$R=\frac {n-1} {\varphi^2} \big (\lambda- 2\varphi \varphi''\big ) -(n-1)(n-2)\left ( \frac{\varphi'}{\varphi}\right )^2, \eqno (3.7)$$ where $\lambda$ is the Einstein constant of $(\Sigma_{c_0}, \, \bar g)$. By Lemma 2.2,  we know that $R\ge 0$. If $\lambda\le 0$, then it follows from (3.7) that $\varphi''\leq 0$. So $\varphi$ is a positive weakly concave function on $\mathbb R$, hence must be a constant function. Therefore, $$(M^n, g)=({\mathbb R}, \,dr^2)\times\,(\Sigma_{c_0}, \, \bar g).$$ 
 Furthermore, it is necessary that $\lambda=0$ since $(M^n, g)$ is a steady Ricci soliton. 

Next, we rule out the possibility of Einstein constant $\lambda>0$.  For any gradient steady soliton with the warped metric of the form (3.1), by (3.4)-(3.6), the gradient steady Ricci soliton equation (1.1) reduces to the following ODE system for the potential function $F(r)$ and the warping function $\varphi (r)$:  
\[ F''=-(n-1)\frac {\varphi''} \varphi,  \qquad \varphi\varphi'F'=\lambda-\varphi \varphi''-(n-2)(\varphi')^2, \eqno (3.8)\] which can be further reduced to a second order nonlinear ODE in $\varphi'$.  By (3.8), (3.7) and Lemma 2.2, it is easy to observe that $\varphi'$  cannot vanish anywhere. Thus, $\varphi'$ does not change sign and the warping function $\varphi$ is monotone over ${\mathbb R}$. Without loss of generality, we assume $\varphi$ is monotone increasing in 
$r$.  However, the fact that $\varphi$ is monotone increasing in $r$ and $\varphi>0$ on ${\mathbb R}$ implies $\lim\inf\varphi'(r)= 0$ as $r\to -\infty$, again a contradiction to (3.8), (3.7) and Lemma 2.2. 
Alternatively, if it were true that the Einstein constant $\lambda>0$ and that $(M^n, g_{ij})=({\mathbb R}, \,dr^2)\times\!{_\varphi}(\Sigma_{c_0}, \, \bar g) $ is a complete gradient steady soliton, then by replacing $(\Sigma_{c_0}, \, \bar g)$ by $({\mathbb S}^{n-1}, \, \bar g_0)$ 
one would produce a complete rotationally symmetric gradient steady soliton ${\mathbb R}\times\!{_\varphi}({\mathbb S}^{n-1},\bar g_0)$ with the same potential function $F$. But this is a contradiction to the classification result of \cite{CC09, CM11} for locally conformally flat gradient steady solitons. Therefore, $\lambda>0$ cannot occur. 

\medskip
In conclusion, $(M^n, g_{ij})$ is either Ricci flat, or isometric to the Bryant soliton, or isometric to the product steady soliton $N^{n-1} \times {\mathbb R}$, where $N^{n-1}$ is Ricci-flat. This completes the proof of Theorem 1.1. 

\qed

\medskip

Clearly,  similar arguments as in the proof of Theorem 1.1 can be used to treat gradient shrinking and expanding Ricci solitons with vanishing $D$-tensor.  

In the shrinking case (i.e., $\rho>0$), by using the same method as in the proof of Theorem 1.1, we have an alternative proof of the following result, which is also a consequence of Proposition 2.3 and the rigidity theorem for gradient shrinking Ricci solitons with harmonic Weyl tensor due to Fern\'andez-L\'opez and Garc\'ia-R\'io \cite{FG1}, and Munteanu-Sesum \cite{MS}. 

\begin{proposition} Let $(M^n, g_{ij}, f)$, $n\ge 5$, be a complete 
gradient shrinking Ricci soliton with vanishing $D$-tensor. Then,  $(M^n, g_{ij}, f)$ is either Einstein, or a finite quotient of the Gaussian soliton ${\mathbb R}^n$, or a finite quotient of the product shrinking soliton $N^{n-1} \times {\mathbb R}$, where $N^{n-1}$ is an Einstein manifold with positive Ricci curvature. 
\end{proposition}

In the expanding case (i.e., $\rho<0$), we have

\begin{proposition}  Let $(M^n, g_{ij}, f)$, $n\ge 5$, be a complete noncompact 
gradient expanding Ricci soliton with vanishing $D$-tensor. Then,  $(M^n, g_{ij}, f)$ is either

\smallskip
(a) Einstein (of negative scalar curvature) with a constant potential function; or

\smallskip
(b) rotationally symmetric, hence  a quotient of an expanding soliton of the form
$$\large([0, \infty\large), \,dr^2)\times\!{_\varphi} ({\mathbb S}^{n-1},\bar g_{0}),$$ where $\bar g_{0}$ is the round metric 
on ${\mathbb S}^{n-1}$; or

\smallskip
(c) a quotient of some warped product expanding Ricci soliton of the form
$$({\mathbb R, \,dr^2})\times\!{_\varphi} (N^{n-1}, {\bar g}),$$ where $(N^{n-1}, {\bar g})$ is an Einstein manifold of negative scalar curvature.  
\end{proposition}

\medskip
\noindent 
{\bf Corollary 3.4.} {\em Let $(M^n, g_{ij}, f)$, $n\ge 5$, be a complete noncompact 
gradient expanding Ricci soliton with nonnegative scalar curvature $R\ge 0$ and vanishing $D$-tensor. Then $(M^n, g_{ij}, f)$ is rotationally symmetric and isometric to a quotient of an expanding soliton of the form
$$\large([0, \infty\large), \,dr^2)\times\!{_\varphi} ({\mathbb S}^{n-1},\bar g_{0}).$$

\end{document}